\documentclass[11pt, a4paper]{article}

\usepackage[utf8]{inputenc}
\usepackage{amsmath, amsthm, amsfonts, amssymb}

\usepackage[bookmarks]{hyperref}
\usepackage{indentfirst} 

\usepackage{authblk}
\usepackage{lipsum}

\providecommand{\keywords}[1]
{
  \small	
  \textbf{\textit{Keywords---}} #1
}

\newtheorem{result}{Theorem}

\theoremstyle{definition}

\newcommand{\abs}[1]{\ensuremath{\left| #1 \right|}}
\newcommand{\op}{\operatorname}

\newcommand{\ze}[1]{\operatorname{Z}(#1)}

\newcommand{\rad}[2]{\op{O}_{#1}(#2)}

\newcommand{\syl}[2]{\op{Syl}_{#1}\left(#2\right)}
\newcommand{\hall}[2]{\op{Hall}_{#1}\left(#2\right)}

\newcommand{\groupgen}[1]{\langle #1 \rangle}

\newcommand{\irr}{\operatorname{Irr}}

\newcommand{\gal}{\operatorname{Gal}}

\newcommand{\Bcd}[1]{\op{Bcd}}

\renewcommand{\phi}{\varphi}
\renewcommand{\epsilon}{\varepsilon}

\title{A note on cut groups of odd order}

\date{\today}

\author{Nicola Grittini}

\affil{Università degli Studi di Firenze}

\begin{document}

\maketitle

\begin{abstract}
A group $G$ is said to be cut if, for every $g \in G$, each generator of $\groupgen{g}$ is conjugated to either $g$ or $g^{-1}$. It is conjectured that a Sylow 3-subgroup $P$ of a cut group $G$ is cut. We prove that this is true if $\abs{G}$ is odd.

\end{abstract}\hspace{10pt}

\section{Introduction}

The concept of cut groups arises form the study of group rings. A finite group $G$ is said to be a \textit{cut group} if $\mathbb{Z}G$ only contains trivial central units (see \cite[Definition~1.1]{Bachle:Integral_group_rings} for details). There exists, however, also some purely group theoretical conditions for a group to be cut. By \cite[Proposition~2.2]{Bachle:Integral_group_rings}, $G$ is a cut group if and only if, for every element $g \in G$, each generator of the cyclic group $\groupgen{g}$ is conjugated either to $g$ or to $g^{-1}$. Equivalently, a group is cut if, for each $\chi \in \irr(G)$, $\mathbb{Q}(\chi) = \mathbb{Q}(\sqrt{-d})$ for some non-negative integer $d$, where $\mathbb{Q}(\chi)$ is the field of values of $\chi$. In particular, for every $\chi \in \irr(G)$, $\abs{\mathbb{Q}(\chi) : \mathbb{Q}} \leq 2$.

Notice that, in the literature, cut groups are also called \emph{inverse semi-rational groups} (see, for instance, \cite{Chillag-Dolfi:Semi-rational}).

It is conjectured that, if $G$ is a cut group, then a Sylow 3-subgroup $P$ of $G$ must also be cut. The conjecture have already been studied for several classes of groups. In particular, in \cite[Theorem 6.6]{Bachle-Caicedo-Jespers-Maheshwary:Global_local_properties} it is proved that it holds when the group is supersolvable, Frobenius or simple. It is also proved that a Sylow 3-subgroup of a cut group $G$ is cut when $\abs{G}$ is odd and $\rad{3}{G}$ is abelian. We will prove that the property holds also without assumptions on $\rad{3}{G}$.

\begin{result}
\label{result:odd_order}
Let $G$ be a cut group of odd order and let $P$ be a Sylow 3-subgroup of $G$. Then, $P$ is cut. Moreover, also $\rad{3}{G}$ is cut.
\end{result}

\section{Proof of Theorem~\ref{result:odd_order}}

We now prove Theorem~\ref{result:odd_order}. The reader shall keep in mind that a character $\phi \in \irr(P)$ of a $p$-group $P$ has values in $\mathbb{Q}_{p^a}$, i.e., the $p^a$-cyclotomic extension of $\mathbb{Q}$, for some $a \in \mathbb{N}$ such that $p^a=\abs{P}$. Therefore, in order for a Sylow 3-subgroup $P$ to be cut, we need that $\phi$ has values in $\mathbb{Q}_3 = \mathbb{Q}(\sqrt{-3})$, since it is the only subfield of $\mathbb{Q}_{3^a}$ to be an extension of $\mathbb{Q}$ of degree 2, thus, the only one which can be generated by the square root of a negative integer.

\begin{proof}[Proof or Theorem~\ref{result:odd_order}]
Let $G$ be a cut group of odd order; it follows from \cite[Remark 13]{Chillag-Dolfi:Semi-rational} and \cite[Theorem 3]{Chillag-Dolfi:Semi-rational} that $G$ is either a 3-group, a Frobenius group of order $3 \cdot 7^a$ or a group of order $7 \cdot 3^b$. In the first case there is nothing to prove while, in the second case, the thesis follows from \cite[Theorem 6.6, (2)]{Bachle-Caicedo-Jespers-Maheshwary:Global_local_properties} and from the fact either $\rad{3}{G}=P$ or $\rad{3}{G}=1$.

Therefore, we only have to consider the case when $\abs{G}=7 \cdot 3^b$. In this situation, let $O=\rad{3}{G}$ and let $H \in \hall{3'}{G}$; it follows from \cite[Theorem 3]{Chillag-Dolfi:Semi-rational} that $OH$ is a Frobenius group with Frobenius kernel $O$, $OT \in \syl{3}{G}$ for some group $T$ of order 3 and $G/O$ is the nonabelian group of order 21 (notice that $G/O$ is a Frobenius group, too, and for this reason groups like $G$ are sometimes called 2-Frobenius or double Frobenius groups).

Since $O$ is a $3$-group, $\ze{O} > 1$ and, since $O$ is normal in $G$, then $\ze{O} \lhd G$. It follows that there exists $M \leq \ze{O}$ minimal normal subgroup of $G$. In particular, $M$ is elementary abelian.

Now, let $\phi \in \irr(P)$. If $M \leq \ker \phi$, then $\phi$ is a character of $P/M \in \syl{3}{G/M}$ and it has values in $\mathbb{Q}_3$ by induction on $\abs{G}$. Thus, we may assume that $M \nleq \ker \phi$.

In order to prove that $\phi$ has values in $\mathbb{Q}_3$, we need to prove that it is fixed by every element of $\gal(\mathbb{Q}_{3^b} \mid \mathbb{Q}_3)$. Thus, let $\sigma \in \gal(\mathbb{Q}_{3^b} \mid \mathbb{Q}_3)$ and let $1_M \neq \lambda \in \irr(M)$ be a constituent of $\phi_M$. Notice that $\lambda$ is linear of order 3 because $M$ is elementary abelian; thus, $\lambda$ has values in $\mathbb{Q}_3$ and is fixed by $\sigma$.

Since $OH$ is Frobenius, no elements of $H$ fix any nontrivial element of $M \leq O$. Thus, also $H \cap \op{I}_G(\lambda) = 1$. Moreover, $O \leq \op{I}_G(\lambda)$ because $M$ is central in $O$. Thus, either $\op{I}_G(\lambda) = O$ or $\op{I}_G(\lambda) = P$.

If $\op{I}_G(\lambda) = P$, then $\phi^G = \chi \in \irr(G)$ by Clifford theorem. Moreover, $\chi$ is fixed by $\sigma$, since $\abs{\mathbb{Q}(\chi) : \mathbb{Q}}=2$ by hypothesis and $o(\sigma)$ is odd. Since both $\chi$ and $\lambda$ are fixed by $\sigma$, it follows from the uniqueness in Clifford theory that also $\phi$ is fixed.

Suppose then that $\op{I}_G(\lambda) = O$, let $\theta \in \irr(O)$ be an irreducible constituent of $\phi_O$ lying over $\lambda$ and notice that $\theta^G = \chi \in \irr(G)$ and $\theta^P = \phi$. By the same argument as the previous paragraph, we have that $\theta$ is fixed by $\sigma$ and, since $\sigma$ commutes with the conjugation by elements of $G$, we have that
$$\phi^{\sigma} = (\theta^{\sigma})^P = \theta^P = \phi.$$

It only remains to prove that $O$ is cut. Suppose that there exists $\theta \in \irr(O)$ and $\sigma \in \gal(\mathbb{Q}_{3^b} \mid \mathbb{Q}_3)$ such that $\theta^{\sigma} \neq \theta$ and let $\phi \in \irr(P)$ lying over $\theta$. Let $M \leq O$ be a minimal normal subgroup of $G$ central in $O$, as above, and let $\lambda$ be an irreducible constituent of $\theta_M$. We can assume that $\lambda \neq 1_M$, since otherwise $\theta$ is an irreducible character of $O/M = \rad{3}{G/M}$ and the thesis follows by induction on $\abs{G}$. Thus, we have that no nontrivial element of $H$ fixes $\lambda$, because $OH$ is Frobenius and $M$ is abelian, as above; since $\theta_M=\theta(1)\lambda$, the same is true for $\theta$. Thus, $H \cap \op{I}_G(\theta) = 1$.

Moreover, let $T=\groupgen{t}$ for some $t \in T$ of order 3, so that $P=O\groupgen{t}$. Since $\sigma$ fixes $\phi$, because $P$ is cut, and it does not fix $\theta$, we have that $\phi_O \neq \theta$ and it follows that $\phi_O= \theta + \theta^t + \theta^{-t}$ and, thus, $\op{I}_G(\theta) = O$. However, also $\phi_O= \theta + \theta^{\sigma} + \theta^{\sigma^{-1}}$, since $\theta^{\sigma},\theta^{\sigma^{- 1}}$ both lie under $\phi$. Thus, either $\theta^{\sigma}=\theta^t$ or $\theta^{\sigma}=\theta^{t^{-1}}$ and, without loss of generality, we may assume $\theta^{\sigma}=\theta^t$.

Let $1\neq h \in H$ and notice that $\theta^h$ is not $\sigma$-invariant and, thus, it does not extend to $P$. Therefore, $\phi_1 = (\theta^h)^P$ is irreducible and, for the same arguments as above, it follows that either $(\theta^h)^{\sigma}=(\theta^h)^t$ or $(\theta^h)^{\sigma}=(\theta^h)^{t^{-1}}$. However, we also have that $(\theta^h)^{\sigma} = (\theta^{\sigma})^h = \theta^{th}$. Since $HT$ is a complement for $O = \op{I}_G(\theta)$ in $G$, $HT \cap \op{I}_G(\theta) = 1$ and we have that either $th=ht$ or $th = ht^{-1}$. In the first case we have that $t$ and $h$ commute, while in the second case we have that $t^h = t^{-1}$ and, thus, $H=\groupgen{h}$ normalizes $T$. Since $HT$ is a Frobenius group, both results are absurd; thus, it follows that $\theta$ is fixed by $\sigma$.
\end{proof}


\end{document}